\documentclass[a4paper,12pt]{article}

\usepackage[italian]{babel}
\usepackage[latin1]{inputenc}

\usepackage{amsmath,amsfonts,amssymb}

\newcommand{\R}{\ensuremath{\mathbb R}}

\newcommand{\N}{\ensuremath{\mathbb N}}

\title{How much is convenient to defect ? \\ A method to estimate the cooperation probability in Prisoner's Dilemma and other games}

\date{1st version: 28.8.2010, last update: 19.1.2014}

\author{Cesco Reale  \thanks{Italian Festival of Mathematical Games. Website: www.cescoreale.com/matematica Email: namesurname@gmail.com, insert "cesco" instead of "name" and "reale" instead of "surname".} \\}
\begin{document}

\maketitle

\begin{abstract}
In many cases the Nash equilibria are not predictive of the experimental players' behaviour.
For some games of Game Theory it is proposed here a method to estimate the probabilities with which the different options will be actually chosen by balanced players, i.e. players that are neither too competitive, nor too cooperative. This will allow to measure the intrinsec cooperativeness degree of a game, only in function of its payoffs.
The method is shaped on the Prisoner's Dilemma, then generalized for asymmetric tables, N players and N options. It is adapted to other conditions like Chicken Game, Battle of the Sexes, Stag Hunt and Translators (a new name proposed for a particular condition). Then the method is applied to other games like Diner's Dilemma, Public Goods Game, Traveler's Dilemma and War of Attrition.
These games are so analyzed in a probabilistic way that is consistent to what we could expect intuitively, overcoming some known paradoxes of the Game Theory.\\

En multaj kazoj la Nash ekvilibroj ne anta\u{u}diras la eksperimentan konduton de la ludantoj.
Por iuj ludoj de Ludteorio estas proponata \^ci-tie metodo por taksi la probablojn la\u{u} kiuj la malsamaj opcioj estos fakte elektataj far ekvilibraj ludantoj, t.e. ludantoj kiuj estas nek tro konkuremaj, nek tro kunlaboremaj.
Tio permesos anka\u{u} mezuri la indekson de kunlaboremeco de iu ludtabelo, nur funkcie de \^giaj valoroj.
La metodo estas studita pri la Dilemo de la Kaptitoj, poste \^generaligita pri nesimetriaj tabeloj, N ludantoj kaj N opcioj. Estas adaptita al aliaj kondi\^coj kiel la Ludo de la Kokido, la Batalo de la Seksoj, la \^Caso de la Cervo kaj la Tradukantoj (nova nomo proponita por iu speciala kondi\^co). Poste la metodo estas aplikita al aliaj ludoj kiel la Dilemo de la Vesperman\^go, la Ludo de la Publikaj A\^joj, la Dilemo de la Voja\^ganto kaj la Milito per Eluzi\^go. 
Tiel \^ci-tiuj ludoj estas analizitaj la\u{u} probableca maniero kiu estas kohera kun tio kion oni povus atendi intuicie, preterpasante iujn konatajn paradoksojn de Ludteorio.

\end{abstract}

\section{Nash equilibria are not always predictive}

In many cases the Nash equilibria are not predictive of the experimental players' behaviour.

For instance, ``in the Public Goods Game repeatedly played, experimental observations show that individuals do not play the predicted noncooperative equilibria'', at least not immediately (Ahn \& Janseen, 2003, Adaptation vs. Anticipation in Public-Good Games).

``In the Traveler's Dilemma it seems very unlikely that any two individuals, no matter how rational they are and how certain they are about each other's rationality, each other's knowledge of each other's rationality, and so on, will play the Nash equilibrium'' (Kaushik Basu, "The Traveler's Dilemma: Paradoxes of Rationality in Game Theory"; American Economic Review, Vol. 84, No. 2, pages 391-395; May 1994).

This paradox has led some to question the value of game theory in general, while others have suggested that a new kind of reasoning is required to understand how it can be quite rational ultimately to make non-rational choices.
In this sense, Douglas Hofstadter proposed the theory of Superrationality: ``it is assumed that the answer to a symmetric problem will be the same for all the superrational players. The strategy is found by maximizing the payoff to each player, assuming that they all use the same strategy. In the Prisoner's Dilemma two superrational players, each knowing that the other is also a superrational player, will cooperate''
(Wikipedia, ref. Douglas R. Hofstadter, 1985, "Metamagical Themas", Basic Books).

I will try to quantify this concept, trying to associate to each option a probability that estimates how many balanced players will actually choose that option. A balanced player is a player that is neither too competitive, nor too cooperative; in paragraph 6 we will define it more precisely. This will allow to measure the intrinsec cooperativeness
degree of a game, only in function of its payoffs. 
Only the competitive games and not the cooperative ones will be analyzed.

This study was actually born from a practical need, i.e. preparing some balanced bimatrixes (more simply we will call them ``tables'') for a game of mathematics and diplomacy, based on the Prisoner's Dilemma. We will talk about it after presenting the proposed model.

\section{Prisoner's Dilemma, other games and Hofstadter's Superrationality}

 The Prisoner's Dilemma was originally framed by Merrill Flood and Melvin Dresher working at RAND Corporation in 1950. Albert W. Tucker formalized the game with prison sentence payoffs and gave it the "prisoner's dilemma" name (Poundstone, 1992).
A classic example of the prisoner's dilemma is presented as follows.

    Two suspects are arrested by the police. The police have insufficient evidence for a conviction, and, having separated the prisoners, visit each of them to offer the same deal. If one testifies for the prosecution against the other (defects) and the other remains silent (cooperates), the defector goes free and the silent accomplice receives the full 10-year sentence. If both remain silent, both prisoners are sentenced to only six months in jail for a minor charge. If each betrays the other, each receives a five-year sentence. Each prisoner must choose to betray the other or to remain silent and each one is assured that the other would not know about the betrayal before the end of the investigation. Both care much more about their personal freedom than about the welfare of their accomplice. (Source: Wikipedia).

We will consider a table of Prisoner's Dilemma, where instead of jail years to minimize there are money prizes to maximize; given $a,b,c,d \in \R$, we define the table $(a,b,c,d)$ in the following way:
if both players cooperate, both receive $b$; if both defect, both receive $c$; if one defects and the other cooperates, the first receives $a$ and the second $d$.
Often these values are indicated with $T,R,P,S$, but in this document we prefer to call them $a,b,c,d$ for several reasons, among them the fact that we will deal with another $p$ indicating a probability.

In the rest of the document we will analyze several conditions, sometimes they are studied in the literature as a specific game, with a specific name:
\[
\begin{split}
&1.1) \quad a>b>c,d  \\
\end{split}
\]

In the Prisoner's Dilemma $c>d$:

\[
\begin{split}
&1.2) \quad a>b>c>d  \\
\end{split}
\]

while in the Chicken Game $b> d \ge c$ :

\[
\begin{split}
&1.3) \quad a>b>d \ge c  \\
\end{split}
\]

In the Battle of the Sexes:

\[
\begin{split}
&1.4) \quad a>d>c \ge b  \\
\end{split}
\]

In the Stag Hunt:
\[
\begin{split}
&1.5) \quad b>a>c>d  \\
\end{split}
\]

We will also study an anomalous case, that we will call \emph{the Translators} (see later):

\[
\begin{split}
&1.6) \quad a>c \ge b>d  \\
\end{split}
\]

We will start from the following remark: in a table ($a,b,c,d$) like (100, 51, 50, 0) $b-c$ is so small compared with $a-b$ and $c-d$ that each player will probably defect, according to the Nash equilibrium. Instead, in a table like (101, 100, 1, 0), like that one analyzed by Hofstadter for his Theory of Superrationality, the advantages of defection $a-b$ and $c-d$ are so small that is almost not worth the risk to come out with (1, 1) instead (100, 100), and so each player will tend more to cooperation than to defection.

Now, we will analyze the following situation. We have a large number of rational and balanced players that are going to play one time the Prisoner's Dilemma; they are divided in pairs and every pair plays with the same table $(a,b,c,d)$; we will try to find a probability $p$ that could estimate how many players will cooperate, consistently with the previous considerations.
We set $q=1-p$, as the defection probability. So $\forall a,b,c,d \in \R$ respecting the condition 1.2) we would like to find a $0<p<1$, reaching $p=0$ only in the limit case $b \le c$, and $p=1$ only in the limit case: $a \le b$ and $ c \le d$.
This estimation should depend only on the given parameters ($a,b,c,d$), and not on the history of the game, that's why we are not going to consider iterated games, Fictitious Play or Evolutionary Stable Strategies.

\section{Maximin criterion}

We could try a first way applying the maximin criterion, and analyzing the corresponding mixed strategy; we obtain: $bp+dq=ap+cq$, and then:

\[
\begin{split}
&2.1) \quad p=\frac{c-d}{(c-d)-(a-b)}  \\
\end{split}
\]

But this $p$ under the condition 1.2) is not in $[0,1]$, so in this case the 2.1) does not solve the problem seen in par. 2.
In order to have $0 < p < 1$ we must have either:

 $a-b<c-d$, $c>d$ and  $a<b$, and it happens for example in The Stag Hunt ($b>a>c>d$);
 or:

 $a-b>c-d$, $c < d$  and $a > b$, and it happens for example in The Chicken Game ($a>b>d \ge c$).

 Moreover, as we will see later, the 2.1) does not describe well these situations.
 For example, in the Stag Hunt, for $b \rightarrow \infty$  we find that $p \rightarrow 0$, but it should tend to 1.
 And in the Chicken, for $c-d=0$ we find $p=0$, but if $a-b<<b-c$, defecting does not make sense.

\section{Maximization of expected payoff}

We can try a second way: we can study the $p$ that maximizes the expected payoff $ \mu(p) = p^2b+(1-p)^2c+p(1-p)(a+d)$; we want $0 < p < 1$, and we have $ \mu (0)=c$ and $ \mu (1)=b$.
Putting $ \frac{d\mu}{dp}=0$, we obtain

\[
\begin{split}
&2.2) \quad p=\frac{a+d-2c}{2(a-b-c+d)}= \eta  \\
\end{split}
\]

$\mu(\eta)$ is a maximum if $c-d<a-b$, otherwise it is a minimum. If it is a maximum, then $p_{MAX}$ (the maximizing $p$) equals to $\eta$ and we can see that $b>c \Rightarrow \eta>1/2$.
In order to have $ \eta<1 $ we find $a-b>b-d$.
So, if $a-b>b-d \Rightarrow 1/2<\eta<1  \Rightarrow  p_{MAX} = \eta$ and $\mu(p_{MAX})>b$; in the other cases $p_{MAX} =1$ and $\mu(p_{MAX})=b$.
Also in this second way, under the condition 1.2) the $p$ in 2.2) is not always in $[0,1]$.

\section{Maximization of expected payoff, given the opponent cooperation probability}

A third way could be to try to maximize the payoff $\mu_{1}$ of player 1, given the $p_{2}$, the cooperation probability of player 2.

$\mu_{1} = p_{1}p_{2}b + (1-p_{1})(1-p_{2})c + p_{1}(1-p_{2})d +p_{2}(1-p_{1})a$

$\mu_{1}(p_{1})$ is a linear function with domain $0 \le p_{1} \le 1$, having maximum in $p_{1} = 0$ if the  $p_{1}$ coefficient $p_{2}(b+c-a-d)+(d-c)<0$ and maximum in $p_{1} = 1$ if $p_{2}(b+c-a-d)+(d-c)>0$.

If $p_{2}(b+c-a-d)+(d-c)=0$  we have $p_{2}= \frac{c-d}{(c-d)-(a-b)}=\theta$ and in this case the function is flat.

Studying the function $\mu_{1}(p_{1})$ it is possible to see that  $p_{1MAX}=1$ only if $\theta \in [0,1) (\Leftrightarrow c \le d)$ and $p_{2}<\theta$,
otherwise $p_{1MAX}=0$.
This is consistent with the Nash equilibria ($p_{1MAX}=0$) in the Prisoner's Dilemma ($a>b>c>d$) and also with the non-coordination Nash equilibria in the Chicken Game ($a>b>d \ge c$). In the Chicken Game, assuming that we know $p_{2}$, we can switch our $p_{1}$ from 0 to 1 if we see that $p_{2}<\theta$, and vicecersa. If $p_{2}=\theta$ then our $p_{1}$ does not affect $\mu_{1}$.

But still we did not find what we are looking for, because this third approach does not solve the problem seen in par. 2.

\section{The proposed estimation of the cooperation probability}

Let us examine a fourth way. Under the condition 1.2) ($a>b>c>d$), we try to think the possible behaviour of a player X playing against a player Y. We can assume that the probability $p_x$ that s/he will cooperate is proportional to $b-c$  (the benefit received by two cooperating players compared with two defecting players), while the probability $q_x$ that s/he will defect is proportional to $p_y(a-b) + q_y(c-d)$ (the benefits received by player X defecting instead of cooperating, weighted according to the cooperation and defection probabilities of player Y). We define such a player as {\em balanced \/}, meaning that this behaviour weights equally the competition benefit and the cooperation benefit, so we can say that such a player is neither too competitive, nor too cooperative.  

We could start our reasoning giving to $p_y$ an arbitrary initial value $0 \le p_{0} \le 1$ ($q_{0}=1-p_{0}$).
Said $b-c = \phi$ and $p_{0}(a-b) + q_{0}(c-d) = \chi_{0}$, the first estimation of $p_x$ is  $p_{1} = \phi / (\phi+\chi_{0})$.
Now, using this first estimation of $p_x$, we can try to give a first estimation of $p_y$, considering that, consistently with the previous reasoning, $p_y$ is proportional to $\phi$ and $q_y$ is proportional to $p_{1}(a-b) + q_{1}(c-d) = \chi_{1}$, so $p_{2} = \phi / (\phi+\chi_{1})$.
We can continue this procedure giving a second estimation of $p_x$, then a second estimation of $p_y$, and so on.
Said $p(a-b) + q(c-d) = \chi$, this recursive sequence $p_i$, independently from the starting point $p_0$, will tend to: $p = \phi / (\phi+\chi)$.

From there we obtain a second degree equation:

\[
\begin{split}
&3) \quad p^2(a-b-c+d)+p(b-d)+(c-b)=0  \\
\end{split}
\]
with solution

\[
\begin{split}
&4.1) \quad p=\frac{d-b \pm \sqrt {(b-d)^2+4(b-c)(a-b-c+d)}}{2(a-b-c+d)}  \\
\end{split}
\]
with $a-b-c+d \neq 0$.

If $a-b-c+d = 0$ from the 3) we have more simply

\[
\begin{split}
&4.2) \quad p= \frac{b-c}{a-c}  \\
\end{split}
\]

For example, going back to what we have seen in paragraph 2, in the Hofstadter's table (101, 100, 1, 0) the 4.2) gives $ p=99\% $, and in the table (100, 51, 50, 0) the 4.1) gives $ p \approx 1.96\% $; these results are consistent with what we could expect intuitively.
Hence, we can measure the intrinsec cooperativeness degree of a table, only in function of its payoffs. 
Moreover, given a table and its estimated $p*$, we can say that a group of players played quite cooperatively (or competitively), if their cooperation rate was greater (or smaller) than $p*$.
Or, if a certain player played several times this same table against different players (without knowing their history), we can say that this player played quite cooperatively (or competitively) if his/her cooperation rate was greater (or smaller) than $p*$.

\section{The probability is univocal and always defined}

 We can prove that under condition 1.2) ($a>b>c>d$) there is always only one root of the equation 4.1) in $[0,1]$, so in Prisoner's Dilemma the proposed estimation is univocal and always defined.

We call  $t$ and $v$ the roots of the equation. Then:

$\frac{b-d}{a-b-c+d}=-(t+v)$ and $\frac{c-b}{a-b-c+d}=tv$

Given a table ($a,b,c,d$), all other tables obtained from the first maintaing unchanged the differences $a-b$, $b-c$, $c-d$ are equivalent concerning the calculation of $p$. Therefore we can fix $d=0$  and then express $a$ and $c$ as a function of $b,t,v$.

We obtain $\frac{a}{b}=2-\frac{tv+1}{t+v}$ and $\frac{c}{b}=1-\frac{tv}{t+v}$, with $t+v \neq 0$.

From the condition $1.2$, and being $d=0 \Rightarrow b>0$, we have:

$ \frac{a}{b}>1 \Rightarrow 2-\frac{tv+1}{t+v}>1 \Rightarrow \frac{tv+1}{t+v}<1$.

If $t+v>0 \Rightarrow  t+v>tv+1$.

If $t+v<0 \Rightarrow  t+v<tv+1$

and moreover:

$ \frac{c}{b}<1 \Rightarrow 1-\frac{tv}{t+v}<1 \Rightarrow \frac{tv}{t+v}>0$.

If $t+v>0 \Rightarrow  tv>0$

If $t+v<0 \Rightarrow  tv<0$
\\

We first analyze the hypothesis $t+v>0$:

$t+v>0 \Rightarrow tv>0 \Rightarrow t,v>0$.

The roots cannot be both greater than 1. If they were, we would have: $t=1+x, v=1+y$, with $x,y$ real positive.
$t+v>tv+1 \Rightarrow (1+x)+(1+y)>(1+x)(1+y)+1 \Rightarrow  2+x+y>1+x+y+xy+1  \Rightarrow  xy<0$, impossible.

At the same time the roots cannot be both in $[0,1]$:

$t+v>tv+1 \Rightarrow t(1-v)>1-v \Rightarrow$

if $1-v>0 (\Rightarrow v<1) \Rightarrow  t>\frac{1-v}{1-v}=1 $

if $1-v<0 (\Rightarrow v>1) \Rightarrow  t<\frac{1-v}{1-v}=1 $

Hence, in the hypothesis $t+v>0$, there is always one root in $[0,1]$ and the other one greater than 1. \\

We analyze now the other hypothesis $t+v<0$:

$t+v<0 \Rightarrow tv<0$ so the roots have opposite signs.
Said $t<0$ and $v>0$, we have:
$t+v<tv+1 \Rightarrow v(t-1)>t-1 \Rightarrow v<1$.
Therefore in the hypothesis $t+v<0$ there is always one root in $[0,1]$ and the other one negative. \\

Furthermore, it is easy to prove that in the 4.1) the root in $[0,1]$ is always that one with $+\sqrt{(b-d)^2+4(b-c)(a-b-c+d)}$, so:

\[
\begin{split}
&4.3) \quad p=\frac{d-b + \sqrt {(b-d)^2+4(b-c)(a-b-c+d)}}{2(a-b-c+d)}  \\
\end{split}
\]

\section{Equiprobability condition}

In the light of what we said, in some tables it will be more probable that players cooperate and in other it will be more probable that players defect.
We search now under which condition is the boundary between these two types of tables.
So we impose $p=1/2$, therefore

\[
\begin{split}
&5) \quad p=\frac{d-b + \sqrt {(b-d)^2+4(b-c)(a-b-c+d)}}{2(a-b-c+d)}=\frac{1}{2}  \\
\end{split}
\]

We obtain the equiprobability condition:

\[
\begin{split}
&6) \quad a-d=3(b-c) \quad \Leftrightarrow \quad b-c=\frac{(a-b)+(c-d)}{2} \\
\end{split}
\]

so $b-c$ is the average of $a-b$ and $c-d$, as it could be easily deduced also from the initial model, putting $p=q$ (and so $\phi=\chi$).
It is possible to prove that, as expected, $p>1/2 \Leftrightarrow  b-c>\frac{(a-b)+(c-d)}{2}$.

\section{Prisoner's Dilemma with n players}

In the case with 3 players, we define the table as follows:
if everyone cooperates, everyone receives $g$; if one defects and two cooperate, the first one receives $f$, and the other two receive $j$; if two defect and one cooperate, the first two receive $h$ and the third one receives $m$; if everyone defects, everyone receives $k$:

C,C,C  $ \rightarrow g,g,g$;
D,C,C  $ \rightarrow f,j,j$;
D,D,C  $ \rightarrow h,h,m$;
D,D,D  $ \rightarrow k,k,k$;

with $f,g,h,j,k,m \in \R$.

Supposing that one player will cooperate, the table for the other two players becomes: ($f,g,h,j$); studying this table under the condition 1.2) we need $f>g>h>j$.
Supposing that one player will defect, the table for the other two players becomes: ($h,j,k,m$); studying this table under the condition 1.2) we need $h>j>k>m$.
Therefore, the condition 1.2) for 3 players becomes $f>g>h>j>k>m$.

\[
\begin{split}
&7) \quad f>g>h>j>k>m  \\
\end{split}
\]

Now we will analyze the cooperation probability $p$ and the defection probability $q$ consistently to the model for 2 players.
Supposing that one player will cooperate, in the table ($f,g,h,j$)  $p$ is proportional to $g-h=\phi_p$ and  $q$ is proportional to $p(f-g)+q(h-j)=\chi_p$;
supposing that one player will defect,  in the table ($h,j,k,m$)  $p$ is proportional to $j-k=\phi_q$ and $q$ is proportional to $p(h-j)+q(k-m)=\chi_q$.

Hence we consider  $p$ proportional to $p\phi_p+q\phi_q = p(g-h)+q(j-k) = \psi$
and $q$ proportional to  $p\chi_p+q\chi_q =  p^2(f-g) + 2pq(h-j) + q^2(k-m)= \omega$.

Hence $p = \psi / (\psi+\omega)$.

From there we obtain a third degree equation:

\[
\begin{split}
8) \quad  & p^3(f-g-2h+2j+k-m)+p^2(g+h-3j-k+2m)+ \\ & +p(-g+h+2j-k-m)+(-j+k)=0  \\
\end{split}
\]

Putting $p=1/2$, and so  $\psi=\omega$, we obtain the equiprobability condition:

\[
\begin{split}
&8.1) \quad f-m+4(h-j)=3(g-k)  \\
\end{split}
\]

The procedure can be recursively extended to $n$ players, even if for $n\ge 5$, as known, only in some cases we will be able to calculate exactly the roots.

\section{Prisonner's Dilemma with asymmetric tables}

We will now extend the model to asymmetric tables, analogously to paragraph 6.
We define an asymmetric table  $(a_x,b_x,c_x,d_x,a_y,b_y,c_y,d_y)$ in the following way:
if both players cooperate, player X receives $b_x$ and player Y receives $b_y$; if both defect, player X receives $c_x$ and player Y receives $c_y$; if player X defects and player Y cooperates, the first receives $a_x$ and the second $d_y$; if player X cooperates and player Y defects, the first receives $d_x$ and the second $a_y$.

We try to think the possible behaviour of a generical player X playing against a generical player Y. Under the condition 1.2) ($a_x>b_x>c_x>d_x$ and $a_y>b_y>c_y>d_y$), we can assume that the probability $p_x$ that player X will cooperate is proportional to $b_x-c_x$, while the probability $q_x$ that s/he will defect is proportional to $p_y(a_x-b_x) + q_y(c_x-d_x)$. We could start our reasoning giving to $p_y$ an arbitrary initial value $0 \le p_{y0} \le 1$ ($q_{y0}=1-p_{y0}$).

Said $b_x-c_x = \phi_x$ and $p_{y0}(a_x-b_x) + q_{y0}(c_x-d_x) = \chi_{x0}$, the first estimation of $p_x$ is  $p_{x1} = \phi_x / (\phi_x+\chi_{x0})$.
Now, using this first estimation of $p_x$, we can try to give a first estimation of $p_y$, considering that, consistently with the previous reasoning, $p_y$ is proportional to $b_y-c_y =\phi_y$ and $q_y$ is proportional to $p_{x1}(a_y-b_y) + q_{x1}(c_y-d_y) = \chi_{y1}$, so $p_{y2} = \phi_y / (\phi_y+\chi_{y1})$.
We can continue this procedure giving a second estimation of $p_x$, then a second estimation of $p_y$, and so on.
Said $p_{y}(a_x-b_x) + q_{y}(c_x-d_x) = \chi_x$ and $p_{x}(a_y-b_y) + q_{x}(c_y-d_y) = \chi_{y}$, this recursive sequence, independently from the starting point, will tend to the pair ($p_x, p_y$): $p_x = \phi_x / (\phi_x+\chi_x)$, $p_y = \phi_y / (\phi_y+\chi_y)$.

Solving this system we obtain the following second degree equations:

\[
\begin{split}
9.1) \quad & p_x^2(b_x-d_x)(a_y-b_y-c_y+d_y)+ p_x (a_x(b_y-c_y) -a_y(b_x-c_x) +b_x b_y + \\
& +d_x d_y -c_y d_x +c_x d_y -2b_x c_y -2b_y c_x +2b_x d_y ) +(b_x-c_x)(b_y-d_y)=0  \\
9.2) \quad & p_y^2(b_y-d_y)(a_x-b_x-c_x+d_x)+ p_y (a_y(b_x-c_x) -a_x(b_y-c_y) +b_y b_x + \\ & +d_y d_x -c_x d_y +c_y d_x -2b_y c_x -2b_x c_y +2b_y d_x )
+(b_y-c_y)(b_x-d_x)=0  \\
\end{split}
\]

\section{Translators condition}

We come back to the case for 2 players.
We analyze some other conditions that we will find later in some applications and that are useful to define some boundaries of the proposed estimation.

We will examine first the condition 1.6) $a>c \ge b>d$.
We can model this situation in the following way. There are two translators that have an excellent level in Esperanto, Galician and Sardinian; the first one has to translate a text from Esperanto into Galician and the second the same text from Esperanto into Sardinian, but the first is much faster in Galician and the second in Sardinian (same speeds, but swapped). Collaborating means helping the other to translate a certain paragraph of the text. So of course if one is helped by the other he will finish earlier, and with a same excellent level of the translated text; but if they both collaborate, they will need globally more time than if they don't.
So saying that $a,b,c,d$ represent the free time remaining after work in the different cases, we obtain the 1.6).

When the condition 1.2) is not respected, the equations 4.3) and 4.2) give often results not in [0,1]. This happens because, even when $\phi$ and $\chi$ remain positive, some addends are negative. But consistently with the proportional method that we are using, each addend in $\phi$ and $\chi$ must be positive.
In condition 1.6), being $b-c<0$ we need to tune our method considering:  $ \phi =0 $ and $ \chi = c-b + p(a-b) + q(c-d)$, therefore  $p=0$, as we could have expected.

\section{Stag Hunt condition}

An interesting case is the condition 1.5) $b>a \ge c>d$, called Stag Hunt.
In ``Discours sur l'origine de l'inégalité parmi les Hommes'' (1754) Jean-Jacques Rousseau described a situation in which two individuals go out on a hunt: ``Voilà comment les hommes purent insensiblement acquérir quelques idées grossières des engagements mutuels, et de l'avantage de les remplir mais seulement autant que pouvait l'exiger l'intérêt présent et sensible ; car la prévoyance n'était rien pour eux, et, loin de s'occuper d'un avenir éloigné, ils ne songeaient même pas au lendemain. S'agissait-il de prendre un cerf, chacun sentait bien qu'il devait pour cela garder fidèlement son poste ; mais si un lièvre venait à passer à la portée de l'un d'eux, il ne faut pas douter qu'il le poursuivit sans scrupule, et qu'ayant atteint sa proie il ne se soucia fort peu de faire manquer la leur à ses compagnons''.
So, if both collaborate in hunting the stag, they both receive $b$, if both will hunt the less worthy hare both receive $c$ ($c<b$), if one hunts the hare while the other remains alone trying to hunt the stag, the first one receives $a$ ($b>a \ge c$) and the second one receives $d$ ($d<c$).
Here we have $ \phi = b-c + p(b-a) $ and $ \chi = q(c-d)$.
We could expect that under this condition should be always $p=1$.
The equation in this case can be expressed as:

\[
\begin{split}
&10.1) \quad (p-1)(p-\frac{c-b}{-a+b-c+d})=0   \\
\end{split}
\]
with $-a+b-c+d \neq 0$, therefore actually one root is always 1.

Say the second root $r_2= (c-b)/(-a+b-c+d)$. If $r_2 \ge 1$ the attractor of the recursive sequence for the set $[0,1]$ is 1, and so $p=1$.
If $0 \le r_2 <1$ we find

\[
\begin{split}
&10.2) \quad 0 \le \frac {b-c}{a-d} < 1/2   \\
\end{split}
\]

In this case $r_2$ is the attractor for the set  $[0,1]$, and so $p= r_2$.
This result could be unexpected, as $a<b$, but actually it is consistent with the problem: for example, in the condition $a-d>>b-c$ (that respects the 10.2) the risk to receive $d$ is not comparable to the small advantage $b-a$, so $p$ is small.
If $a-b+c-d = 0$, we find $ \chi = q(b-a) \Rightarrow  p=\frac {b-c}{b-c}=1$, and it agrees with 10.2), being $\frac {b-c}{a-d}=1 \Rightarrow r_2 \ge 1 $.

\section{Chicken condition}

In the Chicken Game two drivers run one against the other. If one swerves, he will be considered as a ``chicken'' and the other will win the game. If nobody swerves, they will have a serious crash, and it will be even worse than being considered a chicken. If both swerve, they will share the shame (having a better outcome than being the only chicken).
In this game the condition is 1.3) $a>b>d \ge c$; being $c-d \le 0$ we have  $\phi = b-c +q(d-c)$ and $\chi=p(a-b)$, therefore:

\[
\begin{split}
&10.3) \quad p=\frac{3c-b-2d \pm \sqrt {(3c-b-2d)^2-4(2c-b-d)(a-b+c-d)}}{2(a-b+c-d)}  \\
\end{split}
\]
with $a-b+c-d \neq 0$.

If $a-b+c-d = 0$:
\[
\begin{split}
&10.4) \quad p= \frac{a-c}{2a-b-c}  \\
\end{split}
\]

\section{Battle of the Sexes condition}

A husband and a wife agree to meet this evening, but cannot recall if they will be attending the opera or a boxing match. He prefers the boxing match and she prefers the opera, though both prefer being together to being apart. Let us consider the cooperation probability $p$ as the probability to choose the event preferred by the other.
In this game the condition is the 1.4) ($a>d>c \ge b$). Therefore $\phi=q(d-c)$ and $\chi=c-b+p(a-b)$ and we obtain:

\[
\begin{split}
&10.5) \quad p=\frac{c+b-2d \pm \sqrt {(c+b-2d)^2+4(d-c)(a-b+c-d)}}{2(a-b+c-d)}  \\
\end{split}
\]
with $a-b+c-d \neq 0$.

If $a-b+c-d = 0$ (and so $a-b=d-c$), we have $\phi = q(a-b)$ and then more simply

\[
\begin{split}
&10.6) \quad p= \frac{a-b}{a+d-2b}  \\
\end{split}
\]

We examine here the case $(3,0,0,2)$.

As already seen in equation 2.1, the classical maximin mixed strategy gives $ p_{mm}=\frac{a-c}{(a-b)-(c-d)}$, in this case $ p_{mm}=60\% $, so each player attends their preferred event less often than the other event.
Applying the 10.5) we find $p \approx 45 \%$, so  each player attends their preferred event more often than the other event.

\section{Application to a game based on the Prisoner's Dilemma}

In this game there are 4 to 10 teams (with 1 or more players) representing some nations.
At each round the nations are coupled and play a series of phases with tables of Prisoner's Dilemma (or Chicken), accumulating money that they will be able to invest in weapons or (at the end of the round) in development points. The following parts of the round are dedicated to alliances, spying and wars; the war winners steal money from the losers. After a certain number of rounds, the winner is the nation with most development points.
If the number of teams is odd, at each round there will be a group of 3 nations that will play tables for 3 players.
Hence it was necessary to have in each phase a table for 2 players and a table for 3 players, and in each phase these 2 tables should be equilibrated between them in terms of $p$ and $\mu$ (average payoff in the table, weighted in terms of $p$).
Moreover, having several tables in each round, it was necessary that these tables were equilibrated among them.
The tables were tuned applying the formulas presented in this study; testing the game, the obtained results were very satisfying.
Here we show some of these tables.

In the case for 2 players:

\[
\begin{split}
&11.1) \quad \mu= p^2b+q^2c+pq(a+d) \\
\end{split}
\]

In the case for 3 players:

\[
\begin{split}
&11.2) \quad \mu= p^3g+q^3k+p^2q(f+2j)+pq^2(2h+m)\\
\end{split}
\]

Table 1 : Management of industry and commerce

Example of two companies that don't have the right to agree on a common commercial strategy and that wonder if they should (D) or not (C) lower the prices to conquer parts of the market from their opponent.

1)	D,D = 5,5  	

2) C,C = 7,7	

3) D,C = 10,1	

($p$=35 \%, $\mu$ =5,47)

For 3 players :

1)	D,D,D = 4,4,4	

2) C,C,C = 8,8,8

3) D,D,C = 7,7,2		

4) D,C,C = 10,5,5

($p$=33\%, $\mu$ =5,33) \\

%
%
%
%
%
%
%
%
%
%
%

Table 2 :
Corruption

Every country did small or big illegal actions... Are you going (D) or not (C) to denounce your opponent ?

1) D,D = -2,-2  	

2) C,C = 2,2	

3) D,C = 8,-4 	

($p$=50\%, $\mu$ =1)

For 3 players :

1) D,D,D = -2,-2,-2	

2) C,C,C = 4,4,4		

3) D,D,C = 1,1,-4	

4) D,C,C = 10,-2,-2

($p$=50\%, $\mu$ =1,25) \\

%
%
%
%
%
%
%
%
%
%
%
%
%

Table 5 :
Sharing of subventions

The international funds finance the nations with sums that can be shared more (C) or less (D) honestly.

1) D,D = 5,5  	

2) C,C = 8,8	

3) D,C = 9,2 	

($p$=63\%, $\mu$ =6,43)

For 3 players :

1)	D,D,D = 3,3,3	

2) C,C,C = 8,8,8		

3) D,D,C = 7,7,2		

4) D,C,C = 9,6,6

 ($p$=63\%, $\mu$ =6,63)  \\

\section{Application to the Diner's Dilemma}

In the Diner's Dilemma N friends go to a restaurant and before to order they decide to divide the bill in equal parts.
It is possible to choose between an expensive menu and a cheap menu.
We define $r$ the cost of the expensive menu, $s$ the value of the expensive menu, $u$ the value of the cheap menu, $w$  the cost of the cheap menu, with $r,s,u,w \in \R_+$ and with the condition $r>s>u>w$. If instead of menus, we consider objects, the "value" can be better understood quantifying it as the price at which the object can be resold.
Moreover, we consider the condition:

\[
\begin{split}
&12) \quad s-r/N>u-w/N  \quad \Leftrightarrow \quad \frac{r-w}{s-u}<N\\
\end{split}
\]

i.e. egoistically it is convenient to order the expensive menu.
$(r-w)/(s-u)$ is the ratio between the difference of the costs and the difference of the values, that we will call for short as "costs-benefits ratio" ($R_{cb}$). \\

Remodeling the problem as in the case of the Prisoner's Dilemma, in the case $N=2$ we will have:

$a=s-r/2 -w/2$

$b=u-w$

$c=s-r$

$d=u-r/2 -w/2$

Considering the condition 1.2) for the Prisoner's Dilemma we find

\[
\begin{split}
&13) \quad a>b \quad \Leftrightarrow \quad  R_{cb}<2\\
\end{split}
\]

that is equivalent to 12) for $N=2$, and

\[
\begin{split}
&14) \quad b>c \quad \Leftrightarrow \quad R_{cb}>1\\
\end{split}
\]

From $c>d$ we obtain again the 13).
From the equiprobability condition 6) we find:

\[
\begin{split}
&15) \quad p>1/2 \quad \Leftrightarrow \quad a-d<3(b-c)  \quad \Leftrightarrow \quad  R_{cb}> \frac{4}{3}\\
\end{split}
\]

Being $a-b=c-d$, from 4.2) we obtain $p$ as a function of $R_{cb}$:

\[
\begin{split}
&16) \quad p= \frac{b-c}{a-c} \quad \Leftrightarrow \quad p=2-2/R_{cb}\\
\end{split}
\]

with domain defined by $1<R_{cb}<2$.
In the limit case $R_{cb}=2$ we find $p=1$, and indeed $a=b$ e $c=d$, so there is no convenience in defecting.
In the other limit case $R_{cb}=1$ we find $p=0$, and indeed $b=c$, so there is no convenience in cooperating. \\

In the case $N=3$ we have:

$f=s-r/3-2w/3$

$g=u-w$

$h=s-2r/3-w/3$

$j=u-r/3-2w/3$

$k=s-r$

$m=u-2r/3-w/3$

Considering the condition 7) for the Prisoner's Dilemma with 3 players, we find:

\[
\begin{split}
&17) \quad (f>g \quad OR \quad h>j \quad OR \quad k>m) \quad \Leftrightarrow \quad  R_{cb}<3\\
\end{split}
\]

that is equivalent to 12) for N=3.

\[
\begin{split}
&18) \quad (g>h \quad OR \quad j>k) \quad \Leftrightarrow \quad  R_{cb}>3/2\\
\end{split}
\]

From the equiprobability condition 9) we have:

\[
\begin{split}
&19) \quad p>1/2 \quad \Leftrightarrow \quad    3(g-k)>f-m+4(h-j)  \quad \Leftrightarrow \quad  R_{cb}> 2\\
\end{split}
\]

From equation 8), being equal to 0 the coefficients of $p^3$ e $p^2$, we obtain an equation of first degree that expresses $p$ as a function of $R_{cb}$:

\[
\begin{split}
&20) \quad p= \frac{k-j}{-g+h+2j-k-m} \quad \Leftrightarrow \quad p=2-3/R_{cb}\\
\end{split}
\]

with the domain defined from $1.5<R_{cb}<3$.

Can we suppose 
that for N players the formula is $p=2-N/R_{cb}$ with domain $N/2<R_{cb}<N$ ?

\section{Application to the Public Goods Game}

We find a very similar result in the Public Goods Game.
In the basic game each player has $r$ Euros ($r\in R$) and decides how much s/he wants to put in a common pot. Then the euros in the pot grow by an interest rate of $k>1$ ($k \in R$), and then they are equally redistributed to the players.
We will analyze the simple case of 2 players, where each player can put either $r$ or 0. We find:

$a=r+kr/2$

$b=kr$

$c=r$

$d=kr/2$.

$a>b \Rightarrow k<2 $

$b>c \Rightarrow  k>1$.

$a-b=c-d \Rightarrow   p= \frac{b-c}{a-c}=2-2/k$.

$p>1/2  \Rightarrow  k>4/3$.

Exactly as in the 13-16 of Diner's Dilemma.

Now we can try to see what happens if each player has not only 2 options, but $N+1$ options ($N \in \N$), having the possibility to put in the pot: $0, r/N, 2r/N,..., ir/N, ..., r$.
Said $ir/N$ and $jr/N$ two possible amounts to put in the pot ($i,j \in \N$), said $T_{ij}$ the table considering the two options $i$ or $j$, we define $p_{ij}$ as the cooperation probability in  $T_{ij}$, meaning the probability to put the larger amount between $ir/N$ and $jr/N$. And we define $q_{ij}=1-p_{ij}$. As $p=2-2/k$, $p$ depends only on $k$, and not on $r$,$i$ or $j$; fixed $k_*$, we define the value $2-2/k_*=p_*$, and $q_*=1-p_*$. We define $p_i= p(ir/N)$ the probability to choose the amount $ir/N$ among the $N+1$ options.
In analogy with the paragraph 6, each $p_i$ is proportional to $ U_i=\sum_{j=0}^{i-1} p_{ij} + \sum_{j=i+1}^N q_{ij}=ip_* + (N-i)q_*$.

So, for each $i$ we calculate the proportionality coefficient $U_i$ adding the probabilities to play the amount $ir/N$ in the table $T_{ij}$ for each other $j$.
In fact, for $j<i$, in the table $T_{ij}$, the probability to play the amount $ir/N$ is $p_*$, and for $j>i$ is $q_*$.
We find that $ \sum_{i=0}^N U_i = N(N+1)/2 = W$. Then the $p_i$ are:

\[
\begin{split}
&21) \quad p_i= \frac{U_i}{W} = \frac{ip_* + (N-i)q_*}{N(N+1)/2}\\
\end{split}
\]

For a numerical example, let us set $r=100$, $N=4$ and so we have 5 options: 0, 25, 50, 75 or 100 euros.

For $k=4/3 \Rightarrow p_*=1/2$, we obtain: $p_i=20\%$ each, because we are in the equiprobability condition.

For $k=3/2 \Rightarrow p_*=2/3$, we obtain:

$p_0 = p(0) = 4/30 \approx 13.3\%$,

$p_1 = p(25) = 5/30 \approx 16.7\%$,

$p_2 = p(50) = 6/30 = 20\%$,

$p_3 = p(75) = 7/30 \approx 23.3\%$ and

$p_4 = p(100) = 8/30 \approx 26.7\%$.

The $p_i$ increase with the amount because $k>4/3$.

For $k=6/5 \Rightarrow p_*=1/3$,  we obtain:

$p_0 = p(0) = 8/30 \approx 26.7\%$,

$p_1 = p(25) = 7/30 \approx 23.3\%$,

$p_2 = p(50) = 6/30 = 20\%$,

$p_3 = p(75) = 5/30 \approx 16.7\%$ and

$p_4 = p(100) =4/30 \approx 13.3\%$.

The $p_i$ decrease with the amount because $k<4/3$.

\section{Application to the Traveler's Dilemma}

This game was formulated in 1994 by Kaushik Basu and goes as follows.
An airline loses two suitcases belonging to two different travelers. Both suitcases happen to be identical and contain identical antiques. An airline manager tasked to settle the claims of both travelers explains that the airline is liable for a maximum of \$100 per suitcase, and in order to determine an honest appraised value of the antiques the manager separates both travelers so they can't confer, and asks them to write down the amount of their value at no less than \$2 and no larger than \$100. He also tells them that if both write down the same number, he will treat that number as the true dollar value of both suitcases and reimburse both travelers that amount. However, if one writes down a smaller number than the other, this smaller number will be taken as the true dollar value, and both travelers will receive that amount along with a bonus/malus: \$2 extra will be paid to the traveler who wrote down the lower value and a \$2 deduction will be taken from the person who wrote down the higher amount. The challenge is: what strategy should both travelers follow to decide the value they should write down?

Say $r$ the maximum value, $s$ the minimum value, $t$ the bonus, with $r>s \ge t>0$ ($r,s,t \in \R$). The two players have $N+1$ options: given $v= (r-s)/N$ they can play $s, s+v, s+2v, ..., s+iv, ..., r$, with $i,N \in \N$.

We will try to apply again the considerations in paragraph 6 to the case with the 2 options $s+iv$ and $s+jv$ ($j \in \N$); said $T_{ij}$ the table considering the two options $i$ or $j$, we define $p_{ij}$ as the cooperation probability in  $T_{ij}$, so the probability to play the biggest value between $s+iv$ and $s+jv$. We obtain the following values, with $i>j$:

$a=s+jv+t$

$b=s+iv$

$c=s+jv $

$d=s+jv-t$.

We find $a>b \Rightarrow i-j<t/v $, and $t>0 \Rightarrow  c>d$.

$b>c \Rightarrow i>j$, already known.

Applying the 4.3, we obtain the cooperation probability $p_{ij}$ (for $i>j$ it is the probability to play $s+iv$):

\[
\begin{split}
&22) \quad p_{ij}=\frac{-(t+ (i-j) v) + \sqrt {(t+(i-j)v)^2-4(i-j)^2 v^2}}{2(j-i)v}  \\
\end{split}
\]

If $i-j \ge t/v \Leftrightarrow b>a $ we are under the condition $b>a>c>d$, so we must apply the 10.1).

If $i<j$ we just swap $i$ and $j$, obtaining the same cooperation probability, that in this case will be the probability to play $s+jv$.
We can see that $p_{ij}$ depends on $\mid i-j \mid$, but not on $i$ or $j$ separately.

From the equation 6, we can see that the equiprobability condition is:

\[
\begin{split}
&23) \quad 3(b-c)>(a-d) \Leftrightarrow  i-j>2t/3v \\
\end{split}
\]

With the same method used for the Public Goods Game, we have $ U_i=\sum_{j=0}^{i-1} p_{ij} + \sum_{j=i+1}^N q_{ij}$, $W=\sum_{i=0}^N U_i$ and $p_i=U_i/W$.

For a simple numerical example, let us set $r=4$, $s=2$, $t=2$, $N=2$ (3 options), $v=1$.
We can see that for $\mid i-j \mid =1$ we have $p_{ij}\approx 38\%$;  for $\mid i-j \mid \ge 2$, considering the 10.2), we can check that $(b-c)/(a-d) \ge 1/2$, hence we have always $p_{ij}=1$. We obtain:

$p_0 = p(2) \approx 20.6\%$,

$p_1 = p(3) \approx 33.3\%$,

$p_2 = p(4) \approx 46.1\%$.

In the original problem, with $r=100$, $s=2$, $t=2$, $N=98$ (99 options), $v=1$; here also for $\mid i-j \mid \ge 2$, $(b-c)/(a-d) \ge 1/2$, hence we have always $p_{ij}=1$. We obtain

$p_{98} = p(100) \approx 2.01 \%$

$p_0 = p (2) \approx 0.0128 \%$

$p_i = p(i+2) \approx i 0.0206 \%$.

These results are consistent with what we could expect intuitively.

In the article  {\em The Traveler's Dilemma\/}  (Basu, Kaushik. Scientific American Magazine; June 2007) experimental results are reported, where  $r=200$, $s=80$, $N=120$ (121 options), $v=1$. For $t=5$ the average amount proposed by the players was $\mu=180$, and for $t=80$ it was $\mu=120$.
Both results are quite far from the Nash equilibrium ($s=80$).
With our method we obtain: for $t=5$, $\mu= \sum_{i=0}^N (s+iv)p_i \approx 160$ and for $t=80$, $\mu \approx 144$.
Our model is not too far from the experimental results.
We can say that for $t=5$ the players played quite cooperatively (because the average amount 180 was greater than the 160 estimated for balanced players) and for $t=80$ the players played quite competitively (because the average amount 120 was smaller than the 144 estimated for balanced players).

\section{Application to the War of Attrition}

The same method used for the Traveler's Dilemma will be now applied to the War of Attrition, with similar results.
This problem was originally formulated by John Maynard Smith in ``Theory of games and the evolution of animal contests'' (1974, Journal of Theoretical Biology 47: 209-221).
The game works as follows: there are 2 players, each makes a bid; the one who bids the highest wins a resource of value $x \in \R_+$. Each player pays the lowest bid. If each player bids the same amount, they will win $x/2$ each.

We will analyze the case with the 2 bid options $i$ and $j$ ($i,j \in \N$), considering $p_{ij}$ as the cooperation probability (so the probability to play the smallest value between $i$ and $j$). We obtain for the table $T_{ij}$ the following values, with $i>j$:

$a=x-j$

$b=x/2-j$

$c=x/2-i $

$d=-j$.

We find $a>b \Rightarrow x >0 $, already known. $b>c \Rightarrow i>j$, already known.

Applying the 4.3, we obtain the cooperation probability $p_{ij}$ (for $i>j$ it is the probability to play $j$):

\[
\begin{split}
&24) \quad p_{ij}=\frac{-x/2 + \sqrt {x^2/4+4(i-j)^2}}{2(i-j)}  \\
\end{split}
\]

If $i<j$ we just swap $i$ and $j$, obtaining the same cooperation probability, that in this case will be the probability to play $i$.
We can see that $p_{ij}$ depends on $\mid i-j \mid$, but not on $i$ or $j$ separately.
Furthermore, differently from the Traveler's Dilemma, $0<p_{ij}<1$ always. For $i-j$ tending to $+\infty$, $p$ tends to 1.

From the equation 6, we can see that the equiprobability condition is:

\[
\begin{split}
&25) \quad 3(b-c)>a-d \Leftrightarrow  i-j>x/3 \\
\end{split}
\]

With the same method used for the Public Goods Game and for the Traveler's Dilemma,
we have $ U_i=\sum_{j=0}^{i-1} q_{ij} + \sum_{j=i+1}^N p_{ij}$, $W=\sum_{i=0}^N U_i$ and $p_i=U_i/W$.

For a simple numerical example, let us set $x=2$, $0 \le i,j \le 4$. We obtain:

$p_0 = p(0) \approx 29.2\%$

$p_1 = p(1) \approx 25.5\%$

$p_2 = p(2) = 20\%$

$p_3 = p(3) \approx 14.5\%$

$p_4 = p(4) \approx 10.8\%$.

Also in this case the results are consistent with what we could expect intuitively.

\section{Comparison with other experimental results}

Beside the paper of (Basu, 2007), two recent papers (Darai and Grätz, 2010; Khadjavi and Lange, 2013) show results that are not far by our estimation.

In the first one a Prisoner's Dilemma with values (J, J/2, 0, 0) is played.
``The dilemma game is played as follows: each player is assigned two balls, one with the word
{\em steal\/} and one with the word {\em split\/} inside. Then both players choose one of the balls and open
them simultaneously. If both players chose the split ball, the jackpot (J) is divided equally
between the two players. If one player chooses steal and the other chooses split, the former gets
the whole jackpot and the latter receives nothing. If both chose steal, both get nothing''.
They observed an average cooperation rate of 54.5\%, where our method for balanced players gives 61.8\%.
But we must  point out that there is a pre-discussion that could increase the cooperation rate:
``Before the players have to decide which strategy to play, they get some additional time, 
roughly 30 seconds, to discuss with each other what they are going to do''.

In the second one, the  Prisoner's Dilemma table is (9,7,3,1), and in a first group the cooperation rate was 37\%, and in a second group 55\%, where our method for balanced players gives 66\%.

In both papers the results are far from Nash equilibria and much closer to our estimation, and in both cases we can say that the players played more competitively than balanced players.

\section{Conclusions}

The proposed approach seems to describe quite well some classical games of the game theory, using an estimation of the behaviour of balanced players to solve some known paradoxes of the Game Theory.
This estimation can be seen as a convenience index for the different options. 
Moreover, it provides a measure of the intrinsec cooperativeness degree of a game, only in function of its payoffs.

It is possible to apply this approach to many other games, only some applications were showed here.
Another interesting result could be to extend this method to calculate the probability density associated to a continuous range of options; for example, in the Public Goods Game, in the Traveler's Dilemma and in the War of Attrition, each player could choose whatever real number in a fixed range.

\section{References}
\quad \quad 1) Ahn \& Janseen.  {\em Adaptation vs. Anticipation in Public-Good Games\/}. 2003 \\

2) Baron; Durieu; Solal. {\em Algorithme de fictitious play et cycles\/}. Recherches Economiques de Louvain, Louvain Economic Review 69(2), 2003 \\

3) Basu, Kaushik. {\em The Traveler's Dilemma: Paradoxes of Rationality in Game Theory\/}. American Economic Review, Vol. 84, No. 2, pages 391-395; May 1994. \\

4) Basu, Kaushik.  {\em The Traveler's Dilemma\/}.  Scientific American Magazine; June 2007. \\

5) Colman, Andrew M. {\em Prisoner's Dilemma, Chicken, and mixed strategy evolutionary equilibria \/}. Behavioral And Brain Sciences (1995) 18:3 \\

6) Darai; Grätz. {\em Golden Balls: A Prisoner?s Dilemma Experiment  \/}. 2010 \\

7) Diekmann, Andreas. {\em Cooperation in an Asymmetric Volunteer's Dilemma Game Theory and Experimental Evidence \/}. International Journal of Game Theory (1993) 22:75-85 \\

8) Güth, Werner. {\em Experimental Game Theory \/}. \\

9) Hofstadter, Douglas. {\em Metamagical Themas\/}. Basic Books, 1985 \\

10) Khadjavi; Lange.{\em Prisoners and their dilemma\}. 2013 \\  

11) Nash, J. F. (1953). {\em Two-person cooperative games \/}. Econometrica \\

12) Rapoport, A. and A. M. Chammah (1965). {\em Prisoner dilemma: A study in conflict and cooperation\/}. Ann Abor: University of Michigan Press. \\

13) von Neumann, J. and O. Morgenstern (1944). {\em Theory of games and economic behavior\/}. Princeton, NJ.: Princeton University Press. \\

\end{document}